\documentclass{article}
\usepackage{amssymb}
\def\F{{\mathbb F}}

\newtheorem{theorem}{Theorem}
\newtheorem{lemma}{Lemma}
\newtheorem{corollary}{Corollary}

\newcommand{\quash}[1]{}

\begin{document}

\title{Polynomial quotients: Interpolation, value sets and Waring's problem}
\author{Zhixiong Chen\\
School of Mathematics, Putian University, \\ Putian, Fujian
351100, P.R. China\\
ptczx@126.com\\
\\
Arne Winterhof\\
Johann Radon Institute for Computational and Applied Mathematics,\\
Austrian Academy of Sciences, \\
Altenberger Stra\ss e 69, A-4040 Linz,
Austria\\
arne.winterhof@oeaw.ac.at}

\maketitle

\begin{abstract}
For an odd prime $p$ and an integer $w\ge 1$, polynomial quotients $q_{p,w}(u)$ are defined by
$$
q_{p,w}(u)\equiv \frac{u^w-u^{wp}}{p} \bmod p ~~  \mathrm{with}~~ 0
\le q_{p,w}(u) \le p-1, ~~u\ge 0,
$$
which are generalizations of Fermat quotients $q_{p,p-1}(u)$.

First,
we estimate the number of elements $1\le u<N\le p$ for which $f(u)\equiv q_{p,w}(u) \bmod p$ for a given polynomial $f(x)$ over the finite field $\mathbb{F}_p$. In
particular, for the case  $f(x)=x$ we get bounds on the number of fixed points of polynomial quotients.

Second, before we study the problem of estimating the smallest number (called the Waring number) of
summands needed to express each element of $\mathbb{F}_p$  as sum
of values of polynomial quotients, we prove some lower bounds on the size of their value sets,
and then we apply these lower bounds to prove
some bounds on the Waring number using results from bounds on
additive character sums and additive number theory.
\end{abstract}



\section{Introduction}

 For an odd prime $p$ and an integer $u$ with $\gcd(u,p)=1$, the {\it
Fermat quotient $q_p(u)$ \/} is defined as the unique integer
$$
q_p(u) \equiv \frac{u^{p-1} -1}{p} \bmod p ~~  \mathrm{with}~~ 0 \le
q_p(u) \le p-1,
$$
and
$$
q_p(kp) = 0, \qquad k \in \mathbb{Z}.
$$
An equivalent definition  is
\begin{equation}\label{FFF}
q_p(u)\equiv \frac{u^{p-1}-u^{p(p-1)}}{p}\bmod p.
\end{equation}
Many number theoretic and cryptographic questions as well as
measures of pseudorandomness have been studied for Fermat quotients
and their generalizations
\cite{ADS,AW,BFKS,C,Chen,CD,COW,CW,DCH,EM,GW,OS,Sha,S,S2010,S2011,S2011b,SW}.

In particular, for all positive integers $w$, we extend (\ref{FFF})
to define
\begin{equation}\label{poly-def}
q_{p,w}(u)\equiv \frac{u^w-u^{wp}}{p} \bmod p ~~  \mathrm{with}~~ 0
\le q_{p,w}(u) \le p-1, ~~u\ge 0,
\end{equation}
which is called a \emph{polynomial quotient} in \cite{CW-poly}. In
fact $q_{p,p-1}(u)=q_{p}(u)$. We have the following relation between
$q_{p,w}(u)$ and $q_{p}(u)$:
\begin{equation}\label{poly-Fermat:relation}
q_{p,w}(u)\equiv -u^wwq_{p}(u) \bmod p
\end{equation}
for all $u\ge 0$ with $\gcd(u,p)=1$. In particular, we get $q_{p,w}(kp)=0$ if $w\ge 2$ and $q_{p,w}(kp)=k$  if $w=1$. We estimated certain character sums of polynomial
quotients in \cite{CW-poly}. Recently the first author (partly with
other coauthors) also applied polynomial quotients to construct
pseudorandom sequences with good cryptographic properties in
\cite{CG,CNW,DKC}.

In this paper,  first we study interpolation polynomials of polynomial quotients (including the number of fixed points of polynomial quotients) and the size of value sets of polynomial quotients defined in
(\ref{poly-def}). Then we apply results of the size of value sets to study an analogue of the \emph{Waring problem} for
polynomial quotients, that is, the question for the smallest
positive integer $s$, which is called the \emph{Waring number} and
denoted by $g(w,N,p)$, such that the equation
$$
q_{p,w}(u_1)+q_{p,w}(u_2)+\cdots+ q_{p,w}(u_s)\equiv c \bmod p, ~~~
0\le u_1, \ldots,\ u_s<N (\le p),
$$
is solvable for any $c\in\mathbb{F}_p$. If such $s$ does not exist, or equivalently
$q_{p,w}(0)=q_{p,w}(1)=\ldots =q_{p,w}(N-1)=0$,
we put $g(w,N,p)=\infty$. Let $\ell$ be the smallest value with $q_{p,w}(\ell)\not\equiv 0
\bmod p$. Then the Waring number $g(w,N,p)$ always exists if $N>\ell$.
Indeed, it is easy to see that  $g(w,N,p)\le p-1$ for $N>\ell$. For $w=p-1$ (and thus for all $w\not\equiv 0 \bmod p$ by (\ref{poly-Fermat:relation})),
$\ell$ is estimated in \cite{BFKS} by $\ell\le (\log p)^{463/252+o(1)}$ for all $p$, which has recently been improved
to $(\log p)^{7829/4284+o(1)}$ in \cite{Shk}.

Let denote by $F(w,N,p;f(x))$ the number of solutions $0\le u<N$ of $q_{p,w}(u)
\equiv f(u)$ for $f(x)\in\F_p[x]$:
$$
F(w,N,p;f(x))=\#\{u \in \{0, \ldots, N-1\}\ : \ q_{p,w}(u) \equiv f(u)
\bmod p\}, ~~~ N\le p.
$$
In particular, $F(w,N,p;x)$ is the number of fixed points of $q_{p,w}$. We prove upper bounds on $F(w,N,p;f(x))$ in Section \ref{Sec-Interpolation}.

Let denote by
$V(w,N,p)$ the size of the value set of $q_{p,w}(u)$ with $0\le u<N$:
$$
V(w,N,p)=\#\{ q_{p,w}(u)\ :\ u= \ 0, \ldots, N-1\}, ~~~ N\le p.
$$
If $w=kp$ for any positive integer $k$, we have $q_{p,kp}(u)=0$ by
(\ref{poly-Fermat:relation}) and thus $F(kp,N,p;f(x))\le \min\{N, \deg(f(x))\}$, $V(kp,N,p)=1$
and $g(kp,N,p)=\infty$.

For any positive $w$ with $p\nmid w$, write $w=w_1+w_2(p-1)$ with
$1\le w_1\le p-1$ and $w_2\ge 0$. By (\ref{poly-Fermat:relation})
again one can get
$$
q_{p,w_1+w_2(p-1)}(u)\equiv -u^{w_1}(w_1-w_2)q_{p}(u)\equiv
w_1^{-1}(w_1-w_2)q_{p,w_1}(u) \bmod p, ~~~0\le u<p,
$$
and thus for $N\le p$
$$F(w_1+w_2(p-1),N,p;f(x))=F(w_1,N,p;w_1(w_1-w_2)^{-1}f(x)),$$
$$V(w_1+w_2(p-1),N,p)=V(w_1,N,p),$$
and
$$g(w_1+w_2(p-1),N,p)=g(w_1,N,p).$$
 (Note that $w_1\not\equiv w_2
\bmod p$ since  $p\nmid w$.) Hence, we may restrict ourselves to
$1\le w\le p-1$ from now on.

We recall that the classic Waring problem is an important research
field in number theory that investigates the smallest $s$ such
that every element of $\mathcal{R}$ is a sum of $s$ $k$-th powers in
$\mathcal{R}$, where $\mathcal{R}$ is an algebraic structure such as integers,
finite fields, residue rings modulo $m$, polynomial rings, function
fields, etc, see e.g., \cite{LW,W98,W01,WW}. Recently, the second author and
other coauthors considered the Waring problem for Dickson polynomials
in finite fields \cite{GW10,OS11,OTW}.

\section{Interpolation of polynomial quotients}\label{Sec-Interpolation}

In this section we prove bounds on $F(w,N,p;f(x))$. We start with a result which is nontrivial if either $w$ is very large
or $\gcd(w,p-1)$ is moderately large.

\begin{theorem}\label{interpolation}
For $1\le w<p$  and $f(x)\in\F_p[x]$ of degree $d$, let
$$
F(w,N,p;f(x))=\#\{u \in \{0, \ldots, N-1\}\ : \
q_{w,p}(u) \equiv f(u) \bmod p\}, ~~~N\le p.
$$
We have
\begin{eqnarray*}
F(w,N,p;f(x)) & \ll &
 \min\left\{ (p-1-w+d)^{1/4}N^{1/2}p^{1/3}, (p-1-w+d)^{1/8}N^{1/2}p^{3/8}, \right.\\
  & &  \left. \frac{1}{\gcd(w,p-1)}d^{1/4}N^{1/2}p^{4/3}, \frac{1}{\gcd(w,p-1)}d^{1/8}N^{1/2}p^{11/8} \right\}.
\end{eqnarray*}
\end{theorem}
Proof. Using (\ref{poly-Fermat:relation}) we reduce the problem for any $w$ to the case $w=p-1$ (the interpolation of
Fermat quotients), i.e.,
we only need to estimate the number of $0\le u<N$ satisfying
\begin{equation}\label{f-a}
-u^wwq_p(u)\equiv f(u) \bmod p.
\end{equation}
We prove two different bounds.

Bound 1. By (\ref{f-a}) we have $q_p(u)\equiv
-w^{-1}u^{p-1-w}f(u) \bmod p$. Then we get
$$
F(w,N,p;f(x))   \ll \left\{
\begin{array}{ll}
(\deg(x^{p-1-w}f(x)))^{1/4}N^{1/2}p^{1/3}, & 1\le \deg(x^{p-1-w}f(x)) \le p^{1/3}, \\
(\deg(x^{p-1-w}f(x)))^{1/8}N^{1/2}p^{3/8}, & p^{1/3}<\deg(x^{p-1-w}f(x)<p,
\end{array}
\right.
$$
by \cite[Theorem 1]{CW-interpolation}. We remark that the proof of
\cite[Lemma 1]{CW-interpolation} (which deals only with $N=p$) can be easily extended to $N\le p$. The
bound is nontrivial only for $p-w=o(p)$.

Bound 2. The values attained by $u^w\bmod p$ for all $0\le u<p$ are the
same as the values  $u^{\gcd(w,p-1)} \bmod p$. For a fixed
primitive element $\gamma\in\mathbb{F}_p$, we consider the
cyclotomic classes of order $\frac{p-1}{\gcd(w,p-1)}$
\begin{equation}\label{Cj}
C_j=\left\{\gamma^{j+\frac{i(p-1)}{\gcd(w,p-1)}}\bmod p \ : \ 0\le
i< \gcd(w,p-1)\right\}, ~~~ j=0,1,\ldots, \frac{p-1}{\gcd(w,p-1)}-1.
\end{equation}
In fact, the $C_j$'s give a partition of $\mathbb{F}_p^*$. For each
$u\in C_j$, we always have $u^w=\gamma^{jw}$, and the number of
solutions $u\in C_j\cap \{0,\ldots,N-1\}$ of (\ref{f-a}) (hence $q_p(u)\equiv -w^{-1}\gamma^{-jw}f(u)
\bmod p$) is bounded by
$$
\ll \left\{
\begin{array}{ll}
(\deg(f(x)))^{1/4}N^{1/2}p^{1/3}, & 1\le \deg(f(x)) \le p^{1/3}, \\
(\deg(f(x)))^{1/8}N^{1/2}p^{3/8}, & p^{1/3}<\deg(f(x))<p,
\end{array}
\right.
$$
by \cite[Theorem 1]{CW-interpolation} again. So
we have
\begin{eqnarray*}
F(w,N,p;f(x)) & \ll & \frac{p-1}{\gcd(w,p-1)} \min \left\{(\deg(f(x)))^{1/4}N^{1/2}p^{1/3},(\deg(f(x)))^{1/8}N^{1/2}p^{3/8}  \right\} \\
 & \ll & \frac{1}{\gcd(w,p-1)}\min \left\{(\deg(f(x)))^{1/4}N^{1/2}p^{4/3},(\deg(f(x)))^{1/8}N^{1/2}p^{11/8}  \right\}
\end{eqnarray*}
since there are $\frac{p-1}{\gcd(w,p-1)}$ many $C_j$'s. This bound is nontrivial only if $\gcd(w,p-1)\ge p^{5/6}$.  ~\hfill $\Box$

\begin{corollary}\label{fixed-coro}
For $1\le w<p$, the number of fixed points of polynomial quotients
$$
F(w,N,p)=\#\{u \in \{0, \ldots, N-1\}\ : \
q_{w,p}(u) \equiv u \bmod p\}, ~~~ N\le p,
$$
satisfies
$$
F(w,N,p)
 \ll
 \min\left\{ (p-w)^{1/4}N^{1/2}p^{1/3}, (p-w)^{1/8}N^{1/2}p^{3/8}, \frac{N^{1/2}p^{4/3}}{\gcd(w,p-1)}
 \right\}.
$$
\end{corollary}

Besides the cases when $p-w=o(p)$ and $\gcd(w,p-1)\ge p^{5/6}$, there is another nontrivial result if
$\gcd(w-1,p-1)\ge p^{1/2+\varepsilon}$, which includes the important case $w=1$.

\begin{theorem}\label{fixed-thm}
For $1\le w<p$, the number of fixed points of polynomial quotients
$$
F(w,N,p)=\#\{u \in \{0, \ldots, N-1\}\ : \
q_{p,w}(u) \equiv u \bmod p\}, ~~~ N\le p,
$$
satisfies
$$
F(w,N,p) \ll \frac{p^{3/2+\varepsilon}}{\gcd(w-1,p-1)}.
$$
\end{theorem}
Proof. Define
$$
\widetilde{C}_j=\left\{\gamma^{j+\frac{i(p-1)}{\gcd(w-1,p-1)}}\bmod
p \ : \ 0\le i< \gcd(w-1,p-1)\right\}, ~~~ j=0,1,\ldots,
\frac{p-1}{\gcd(w-1,p-1)}-1.
$$
The number of solutions $u\in \widetilde{C}_j\cap
\{0,\ldots,N-1\}$ of $q_p(u)\equiv -w^{-1}u^{-(w-1)}\equiv
-w^{-1}\gamma^{-j(w-1)} \bmod p$ is bounded by $O(p^{1/2+\varepsilon})$ by
\cite[Proposition 2.1]{F}. So we have
$$
F(w,N,p) \ll \frac{p-1}{\gcd(w-1,p-1)}p^{1/2+\varepsilon}\ll
\frac{p^{3/2+\varepsilon}}{\gcd(w-1,p-1)}.
$$
The bound is nontrivial only for $\gcd(w-1,p-1)\gg
p^{1/2+\varepsilon}$ and $N\gg p^{1/2+\varepsilon}$. ~\hfill $\Box$

\section{Size of value sets}
First we prove a bound on $V(p-1,N,p)$, the size of the value set of
Fermat quotients $q_{p}$, see \cite[Theorem 13]{OS} for $N=p$. Then we estimate $V(w,N,p)$ for general $1\le
w\le p-2$ in terms of $V(p-1,N,p)$ by (\ref{poly-Fermat:relation}).

\begin{lemma}\label{VS:w=p-1}
Let $ V(p-1,N,p)=\#\{ q_{p}(u)\ :\ u= \ 0, \ldots, N-1\}$. We have
$$
V(p-1,N,p) \gg \frac{N^2}{p\log^{2} N}, ~~~ N\le p.
$$
\end{lemma}
Proof. For $N<p$, one
can get the desired result the same way as for $N=p$, see the proof of
\cite[Theorem 13]{OS}. For the convenience of the reader, we  sketch the proof here.

Let $Q(N,a)$ be the number of primes $l$ smaller than $N$ with
$q_p(l) = a$. Clearly
$$\sum_{a=0}^{p-1} Q(N,a) = \pi(N-1),$$
where $\pi(x)$ denotes the number of primes $l \le x$. The number of
prime number pairs $(l,r)$ with $0\le l,r \le N-1$ and $q_p(l) = q_p(r)$ is
$ \sum_{a=0}^{p-1} Q(N,a)^2$.

According to the fact that $q_p\ :\ \mathbb{Z}_{p^2}^*\rightarrow
\mathbb{Z}_p$ is a group homomorphism with kernel
$\mathrm{ker}(q_p)$ of size $p-1$, we see that $l/r\in
\mathrm{ker}(q_p)$ for each pair $(l,r)$ above.  Now for each $u\in \mathrm{ker}(q_p)$, there
are $\pi(N-1)$ many  pairs $(l,l)$ such that $1\equiv l/l \bmod
{p^2}$ if $u=1$ and at most one pair $(l,r)$ such that $u\equiv l/r
\bmod {p^2}$ if $u\neq 1$, since otherwise, $u\equiv l_1/r_1 \equiv
l_2/r_2 \bmod {p^2}$ leads to $l_1=r_1, l_2=r_2$ or $l_1=l_2,
r_1=r_2$. So we get
$$
\sum_{a=0}^{p-1} Q(N,a)^2\le
\pi(N-1)+\#\mathrm{ker}(q_p)-1=\pi(N-1)+p-2.
$$
On the other hand, only at most $V(p-1,N,p)$ many $Q(N,a)$ are nonzero
for $0\le a\le p-1$, so by the Cauchy-Schwarz inequality we  have
$$
\left(\sum_{a=0}^{p-1} Q(N,a)\right)^2 \le V(p-1,N,p) \sum_{a=0}^{p-1}
Q(N,a)^2.
$$

Putting everything together, we obtain
$$
V(p-1,N,p) \gg \pi(N-1)^2 p^{-1},
$$
which concludes  the proof. ~\hfill $\Box$

As in Section \ref{Sec-Interpolation} we prove different bounds on $V(w,N,p)$ which are nontrivial
if either $\gcd(w,p-1)$ or $\gcd(w-1,p-1)$ is large enough.

\begin{theorem}\label{VS:large w}
For $1\le w<p$  let $ V(w,N,p)=\#\{ q_{p,w}(u)\ :\ u= \ 0, \ldots,
N-1\}$. We have,
$$
V(w,N,p)
 \gg   \gcd(w,p-1)\left(\frac{N}{p\log N}\right)^2, ~~~ N\le p.
$$
\end{theorem}
Proof. The values assumed by $u^w\bmod p$ for all $0\le u<p$ are the
same as the values $u^{\gcd(w,p-1)} \bmod p$. For a fixed
primitive element $\gamma\in\mathbb{F}_p$, we consider the
cyclotomic classes of order $\frac{p-1}{\gcd(w,p-1)}$
defined by (\ref{Cj}). Let $U$ be
the biggest subset of $\{0, \ldots, N-1\}$ such that $q_{p}(u)\neq
q_{p}(v)$ for any $u\neq v\in U$. It is easy to see that $\#
U=V(p-1,N,p)$. Then for any $u_1, u_2\in (C_j \cap U)$ and any $j$,
using (\ref{poly-Fermat:relation}) we always have
$$
u_1^w\equiv u_2^w\bmod p~~ \mathrm{and}~~ q_{p,w}(u_1)\neq
q_{p,w}(u_2).
$$
By the pigeonhole principle we see that there exists  some $j$ with
$$
C_j \cap U\ge \frac{\# U}{(p-1)/\gcd(w,p-1)}.
$$
 So we have
$$
V(w,N,p)\ge \frac{\# U}{(p-1)/\gcd(w,p-1)}\gg
\frac{\gcd(w,p-1)N^2}{p^2\log^2 N}
$$
by Lemma \ref{VS:w=p-1}.  ~\hfill $\Box$

The bound in Theorem \ref{VS:large w} is trivial if $\gcd(w,p-1)\ll \log^2N$. Below we consider the cases of large $\gcd(w-1,p-1)$ (including $w=1$) and get a nontrivial bound using a different
method.

\begin{theorem}\label{VS:large w-1}
For $1\le w<p$  let $ V(w,N,p)=\#\{ q_{p,w}(u)\ :\ u= \ 0, \ldots,
N-1\}$. We have,
$$
V(w,N,p)
 \gg   \gcd(w-1,p-1)\frac{N^{1/2}}{p^{4/3}}, ~~~ N\le p.
$$
\end{theorem}
Proof. We first prove the case $w=1$ and then reduce the general
case $w>1$ to the case $w=1$. The proof follows 
\cite[Section 2]{CW-interpolation}, which deals with the case $N=p$. Put
$$
M_d=\#\{u \in \{0, \ldots, N-1\}\ : \ q_{p,1}(u) =d\}
$$
for some $d$. We first estimate an upper bound on $M_d$.

For $0\le a< N$ and $1\le b< N$,  suppose that $(a,a+b \bmod N)$ is
a pair of points satisfying
$$
q_{p,1}(a) =  q_{p,1}(a+b \bmod N) = d.
$$
 We note that there
are $M_d(M_d-1)$ such pairs. (Note that $M_d=1$ if no such $b$
exists.) Now we fix any $1\le b<N$ and estimate the number of $a$.
For each pair $(a,b)$, set $c=b$ if $a+b<N$ and $c=b-N$ otherwise.
Hence for given $b$ there are two possible choices of $c$ such that
$(a,a+c)$ satisfy
\begin{equation}\label{pair}
q_{p,1}(a) = q_{p,1}(a+c) = d
\end{equation}
for some $a$. For given $c$ we estimate the number of $a$.

If $(a,a+c)$ is a pair satisfying (\ref{pair}), using
(\ref{poly-Fermat:relation}) and the definition of $q_p(u)$ we have
\begin{eqnarray*}
d=q_{p,1}(a+c)\equiv -(a+c)q_p(a+c)&\equiv& -aq_p(a)-cq_p(c)-c\sum_{i=1}^{p-1} \frac{{p \choose i}}{p} (ac^{-1})^{i}\\
&\equiv& q_{p,1}(a)+q_{p,1}(c) +c\sum_{i=1}^{p-1}
\frac{(-ac^{-1})^{i}}{i} \bmod p
\end{eqnarray*}
and thus
$$q_{p,1}(c) +c\sum_{i=1}^{p-1}
\frac{(-a^{-1}c)^{i}}{i}\equiv 0 \bmod p.$$
Substituting $a\equiv
-cx\bmod p$ for $x\in\mathbb{F}_p$ we get
$$
q_{p,1}(c)c^{-1}+\sum_{i=1}^{p-1} \frac{x^i}{i}\equiv 0\bmod p.
$$
Now by \cite[Lemma 4]{H-B}  the number of $x$ (which is not smaller than
 the number of $a$ since $0\le a<N$) for fixed $c$ is
bounded by $O(p^{2/3})$ and we have
$$
M_d(M_d-1)\ll (N-1) \min \{p^{2/3}, N\},
 $$
and thus $M_d\ll N^{1/2}p^{1/3}$ if $N\gg p^{2/3}$, which implies
that
$$
V(1,N,p) \gg   \frac{N^{1/2}}{p^{1/3}}.
$$
From (\ref{poly-Fermat:relation}) again, we have
$$
q_{p,w}(u)\equiv -u^w w q_{p}(u) \equiv u^{w-1} w q_{p,1}(u)\bmod p,
$$
and hence
$$
V(w,N,p)\ge \frac{V(1,N,p)}{(p-1)/\gcd(w-1,p-1)}
$$
following the proof of Theorem \ref{VS:large w}. ~\hfill $\Box$\\

\noindent \textbf{Remark.} Ostafe and Shparlinski stated the problem of finding a nontrivial lower
bound on $V(1,N,p)$ for $N\le p$ in \cite{OS}. In particular,
Theorem \ref{VS:large w-1} implies
$$
V(1,N,p) \gg  N^{1/2}p^{-1/3},
$$
which is nontrivial for $N\gg p^{2/3}$.

\section{Bounds on the Waring number}\label{waring}

\subsection{Bound derived from additive character sums}

We first present a bound on character sums of polynomial quotients, which is
a special case of \cite[Theorem 3]{CW-poly}. In this subsection, we
will exploit these character sums to estimate the Waring number
$g(w,N,p)$.

\begin{lemma}\label{ExpSum}
Let $q_{p,w}(u)$ be defined by (\ref{poly-def}) with $1\le w<p$. For
any nontrivial additive character $\psi$ of $\mathbb{F}_p$ we have,
$$
\left|\sum_{u=0}^{N-1}\psi(q_{p,w}(u)) \right| \ll
\frac{1}{\gcd(w,p-1)} N^{1/2}p^{11/8}, ~~~ N\le p.
$$
\end{lemma}
As noted in \cite[Theorem 2]{H-B}, the exponent $\varepsilon$ in
\cite[Theorem 3]{CW-poly} can be removed when the modulus $k$ of
(multiplicative) characters equals $p^2$ since the Burgess bound
contains a factor $k^{3/16+\varepsilon}$, see \cite[Theorems 2 and
3]{B}. Lemma \ref{ExpSum} is only nontrivial for $N\ge p^{3/4}$. However, using the precise
Theorem 3 in \cite{CW-poly} we can derive bounds which are nontrivial for $N\ge p^{1/2+o(1)}$.

\begin{theorem}\label{th5}
For $1\le w<p$, we have
$$g(w,N,p) \le s \quad \mbox{if}\quad \gcd(w,p-1)^{s-1}\gg  p^{11s/8+1/4}N^{-s/2-1}\log^2 N,\quad s\ge 3, ~~~ N\le p.$$
\end{theorem}
Proof. Without loss of generality we restrict ourselves to the case
 $g(w,N,p)\ge 3$.

Let $\psi$ be a nontrivial additive character of $\F_p$.  For $s\ge
3$ and $y\in \F_p$, the number $N_s(y)$ of solutions
$(v_1,v_2,u_1,\ldots, u_{s-2})$ of the equation
$$
y\equiv v_1+v_2+q_{p,w}(u_1)+\cdots+ q_{p,w}(u_{s-2})\bmod p,~~~
v_1,v_2\in V(w,N,p), ~0\le u_1, \ldots, u_{s-2}<N,
$$
is
\begin{eqnarray*}
N_s(y)&=&\frac{1}{p} \sum_{a\in\F_p}~\sum_{v_1,v_2\in
V(w,N,p)}~\sum\limits_{\stackrel{0\le u_j<N}{1\le j\le s-2}}
   \psi\left(a\left(v_1+v_2+\sum_{i=1}^{s-2}q_{p,w}(u_i) -y\right)\right )\\
&=&
\frac{V(w,N,p)^2N^{s-2}}{p}\\
&& +\frac{1}{p}\sum_{a\in\F_p^*}\psi(-ay)\sum_{v_1,v_2\in
V(w,N,p)}\psi(a(v_1+v_2))~\sum\limits_{\stackrel{0\le u_j<N}{1\le j\le
s-2}}\psi\left(a\sum_{i=1}^{s-2} q_{p,w}(u_i)\right).
\end{eqnarray*}
By Lemma \ref{ExpSum}, we have
\begin{eqnarray*}
\left|N_s(y)-\frac{V(w,N,p)^2N^{s-2}}{p}\right| & \le &
\frac{1}{p}\sum_{a\in\F_p^*}\left|\sum_{v\in V(w,N,p)}\psi(av)
\right|^2\cdot\left |\sum_{0\le u<N}\psi\left(a q_{p,w}(u)\right)
\right |^{s-2}\\
 & \ll &
\frac{1}{p}\cdot \left(\frac{N^{1/2}p^{11/8}}{\gcd(w,p-1)}\right)
^{s-2}\sum_{a\in\F_p^*}\left|\sum_{v\in
V(w,N,p)}\psi(av) \right|^2\\
 &\le &
\frac{1}{p}\cdot \left(\frac{N^{1/2}p^{11/8}}{\gcd(w,p-1)}\right)
^{s-2}\sum_{a\in\F_p}~\sum_{v_1,v_2\in
V(w,N,p)}\psi(a(v_1-v_2)) \\
 &\le & V(w,N,p)\cdot \left(\frac{N^{1/2}p^{11/8}}{\gcd(w,p-1)}\right)
^{s-2}.
\end{eqnarray*}
The number $N_s(y)$ is positive for all $y\in \F_p$ if
$$
V(w,N,p)>p\cdot
\left(\frac{p^{11/8}}{\gcd(w,p-1)N^{1/2}}\right)^{s-2}
$$
and thus $g(w,N,p)\le s$ under this condition. ~\hfill $\Box$

\textbf{Remark}. It is clear that $g(p-1,p,p)\le 3$, which is the Waring number for Fermat quotients. Theorem \ref{th5}
 is only nontrivial if $\gcd(w,p-1)\gg p^{7/8}$.

\subsection{Bound derived from the Cauchy-Davenport
theorem}

In this subsection we prove a bound on $g(w,N,p)$ based on the
Cauchy-Davenport theorem, see e.g., \cite[Theorem 5.4]{TV}, which is rather moderate but nontrivial if
$\gcd(w,p-1)\gg \log^2p$ or $\gcd(w-1,p-1)\gg p^{5/6}$.

\begin{lemma}\label{Cauchy-Davenport}
(Cauchy-Davenport theorem) Let $A,B$ be nonempty subsets of
$\mathbb{F}_p$. Then
$$\#(A+B)\geq \min \{\#A+\#B-1, p\},$$
where $A+B=\{a+b: a\in A, b\in B\}.$
\end{lemma}

\begin{theorem}\label{bound-1}
For $1\le w<p$, we have
$$
g(w,N,p)\ll \min \left\{ \frac{p^3\log^2p}{N^2\gcd(w,p-1)},
\frac{p^{7/3}}{N^{1/2}\gcd(w-1,p-1)}\right\}, ~~~ N\le p.
$$
\end{theorem}
Proof. For $s\ge 1$ put
 $$
 W_s=\{q_{p,w}(u_1)+q_{p,w}(u_2)+\cdots+ q_{p,w}(u_s)\ : \
0\le u_1, \ldots,\ u_s<N\}.
$$
Since $W_s=W_{s-1}+W_1$ for $s\ge 2$, by Lemma
\ref{Cauchy-Davenport} we have
 $$\#W_s\geq \min \{\#W_{s-1}+\#W_1-1, p\},\quad s\ge 2,$$
 and get by induction
 $$\#W_s\geq \min \{s(\#W_{1}-1)+1, p\},\quad s\ge 1.$$
Hence we get
$$
s\le \left\lceil \frac{p-1}{\#W_1-1}\right\rceil
$$
and then the desired result follows from Theorems
\ref{VS:large w} and \ref{VS:large w-1}, respectively.
~\hfill $\Box$\\

\section{Final remarks}

1. The bounds in this paper are non-trivial if $\gcd(w,p-1)$ or
$\gcd(w-1,p-1)$ is ``large". It is challenging to study general
$w$.

2. The bound in Lemma \ref{ExpSum} does not cover the cases of small
$w$. In particular, it is an interesting problem to estimate the character sums
$$
\sum_{u=0}^{N-1}\psi(q_{p,1}(u)).
$$

3. In \cite{S2011} Shparlinski considered the smallest number
$\Lambda_p$ for Fermat quotients such that
$$
\{q_{p}(u)\ : \ u \in \{1, \ldots, \Lambda_p\} \}=\F_p
$$
by estimating $\Lambda_p\le p^{463/252+o(1)}$. It would be
interesting to extend this result to $q_{p,w}$.

4. In \cite{SW}, Shparlinski and the second author introduced the
\emph{polynomial Fermat quotients} in polynomial rings over finite
fields. Let $\F_q$ be a finite field of prime power order $q=p^r$, for a
fixed irreducible polynomial $P\in\F_q[X]$ of degree $n\ge 2$ and
$A\in\F_q[X]$, the polynomial Fermat quotient is defined by
$$
q_P(A)\equiv \frac{A^{q^n-1}-1}{P} \bmod P, ~~~ \deg(q_P(A))<n, ~~
\mathrm{if}~~ \gcd(A,P)=1,
$$
and $q_P(A)=0$ if $\gcd(A,P)=P$. The properties, such as the number
of fixed points and the image size, of the polynomial Fermat
quotient are investigated in \cite{SW}.

Like the definition of polynomial quotients modulo $p$, one can
define
$$
q_{P,w}(A)\equiv \frac{A^{w}-A^{wq^n}}{P} \bmod P, ~~~
\deg(q_{P,w}(A))<n,
$$
for integers $w\ge 1$. In particular, $-q_{P,1}(A)$ has been
introduced in \cite{SS}. Since $q_{P,1}$ is a linear map with kernel
of dimension $\lceil n/p\rceil$, we have
$$
\#\{A: q_{P,1}(A)=B, \deg(A)<n \}=q^{\lceil n/p\rceil}
$$
for any fixed $B=q_{P,1}(A_0)$ for some $A_0$ and hence
$$
\#\{ q_{P,1}(A): \deg(A)<n \}=q^{n-\lceil n/p\rceil}.
$$
(See also the proof of \cite[Lemma 6]{SW}.)

Here we present some lower bounds on the image size of $q_{P,w}$ for
$w>1$. We only consider the case $p\nmid w$, since otherwise
$q_{P,w}$ is a zero map. Firstly from
$$
q_{P,w}(A)\equiv -wA^{w} q_{P}(A) \bmod P,
$$
we reduce the problem to the image size of $q_{P}$ (see
\cite[Theorem 5]{SW}) and obtain
$$
\#\{ q_{P,w}(A): \deg(A)<n \}\gg \frac{\gcd(w,q^{n}-1)}{qn^2}
$$
by using the proof technique of Theorem \ref{VS:large w}. Secondly
from
\begin{equation}\label{PP:w-1}
q_{P,w}(A)\equiv wA^{w-1} q_{P,1}(A) \bmod P,
\end{equation}
we obtain a lower bound similarly in terms of the image size of
$q_{P,1}$ above:
$$
\#\{ q_{P,w}(A): \deg(A)<n \}\gg \frac{\gcd(w-1,q^{n}-1)}{q^{\lceil
n/p\rceil}}.
$$
Finally  from (\ref{PP:w-1}) again, since there are exactly
$\frac{q^n-1}{\gcd(w-1,q^n-1)}+1$ many different $A^{w-1}$ modulo
$P$ for all $A$ with $\deg(A)<n$, then there exists a $B$ such
that at least
$\left(\frac{q^n-1}{\gcd(w-1,q^n-1)}+1\right)/q^{n-\lceil
n/p\rceil}$ many $A$ satisfy $q_{P,1}(A)=B$ but $A^{w-1}  \bmod
P$ are different for all such $A$. Thus we obtain another lower
bound
$$
\#\{ q_{P,w}(A): \deg(A)<n \}\gg \frac{q^{\lceil
n/p\rceil}}{\gcd(w-1,q^n-1)}.
$$

About the Waring problem for $q_{P,w}$, we can not say anything more. The Cauchy-Davenport theorem is not true for
arbitrary fields in general and we do not have any results of
character sums of $q_{P,w}$, so we can not deal with the  Waring
problem using the methods in Section \ref{waring}. But for $q_{P,1}$
the Waring number does not exist, since $q_{P,1}$ is a linear map
with kernel of dimension $\lceil n/p\rceil$ and hence the image of
$q_{P,1}$ is a proper linear subspace of $\F_q[X]/\langle P\rangle$.
That is, there does exist an element in  $\F_q[X]/\langle P\rangle$
which can not be represented as a sum of $q_{P,1}$.

\section*{Acknowledgements}

Z.X.C. was partially supported  by the National Natural Science
Foundation of China under grant No. 61373140 and the Special Scientific Research Program
in Fujian Province Universities of China  under grant No. JK2013044. A.W. was partially supported  by
the Austrian Science Fund (FWF): Project F5511-N26, which is a part of the Special Research Program ``Quasi-Monte
Carlo Method: Theory and Applications".

Parts of this paper were written during  pleasant mutual visits of the authors to
RICAM, Austrian Academy of Sciences in Linz and Putian University. They wish to thank for the hospitality and financial support. The authors also thank Alina Ostafe for pointing to Waring's problem for Fermat quotients.


\begin{thebibliography}{99}

\bibitem{ADS} T. Agoh, K. Dilcher and L. Skula.
Fermat quotients for composite moduli.  J.  Number Theory 66
(1997) 29--50.

\bibitem{AW} H. Aly  and A. Winterhof. Boolean functions derived from Fermat
quotients. Cryptogr. Commun. 3 (2011) 165--174.

\bibitem{BFKS}
J. Bourgain, K. Ford, S. Konyagin and I. E. Shparlinski. On the
divisibility of Fermat quotients. Michigan Math. J. 59 (2010) 313--328.

\bibitem{B} D. A. Burgess. On character sums and L-series, II. Proc. London
Math. Soc.  13   (1963)  524--536.

\bibitem{C}
M. C. Chang. Short character sums with Fermat quotients. Acta Arith.
152 (2012) 23--38.

\bibitem{Chen}  Z. X. Chen. Trace representation and linear complexity of binary sequences derived from Fermat quotients.  http://arxiv.org/arXiv: 1306.5648, 2013.




\bibitem{CD}  Z. X. Chen  and  X. N. Du. On the linear complexity of binary threshold sequences derived from
Fermat quotients. Des. Codes Cryptogr. 67 (2013) 317--323.



\bibitem{CG}Z. X. Chen  and  D. G\'{o}mez-P\'{e}rez. Linear complexity of
binary sequences derived from polynomial quotients. Sequences and Their Applications-SETA 2012, 181--189, Lecture Notes in Comput. Sci., 7280, Springer, Berlin, 2012.

\bibitem{CHD}  Z. X. Chen, L. Hu  and  X. N. Du. Linear complexity of some binary sequences derived from Fermat
quotients. China Commun  9 (2012) 105--108.


\bibitem{CNW}  Z. X. Chen, Z. H. Niu and C. H. Wu. On the $k$-error linear complexity of binary
sequences derived from polynomial quotients. http://arxiv.org/arXiv: 1307.6626, 2013.


\bibitem{COW}  Z. X. Chen, A. Ostafe and  A. Winterhof. Structure of
pseudorandom numbers derived from Fermat quotients. Arithmetic of Finite Fields-WAIFI 2010, 73--85, Lecture Notes in Comput. Sci., 6087, Springer, Berlin, 2010.




\bibitem{CW-poly}Z. X. Chen  and A.  Winterhof. Additive character sums of
polynomial quotients. Theory and Applications of Finite Fields-Fq10,
67--73, Contemp. Math., 579, Amer. Math. Soc., Providence, RI, 2012.

\bibitem{CW}
Z. X. Chen  and A.  Winterhof. On the distribution of pseudorandom
numbers and vectors derived from Euler-Fermat quotients. Int. J.
Number Theory 8 (2012) 631--641.

\bibitem{CW-interpolation} Z. X. Chen and A. Winterhof. Interpolation of Fermat quotients. SIAM J. Discr. Math. 28 (2014) 1--7.



\bibitem{DCH}
X. N. Du, Z. X. Chen and L. Hu. Linear complexity of binary
sequences derived from Euler quotients with prime-power modulus.
Inform. Process. Lett. 112 (2012) 604--609.

\bibitem{DKC}
X. N. Du, A. Klapper and Z. X. Chen. Linear complexity of
pseudorandom sequences generated by Fermat quotients and their
generalizations. Inform. Process. Lett. 112 (2012) 233--237.

\bibitem{EM} R. Ernvall and T. Mets{\"a}nkyl{\"a}.
On the $p$-divisibility of Fermat quotients. Math. Comp.  66 (1997)
1353--1365.


\bibitem{F} W. L. Fouch\'{e}. On the Kummer-Mirimanoff congruences.  Quart.
J. Math.  37 (1986) 257--261.

\bibitem{GW10}
D. G\'{o}mez-P\'{e}rez and A. Winterhof. Waring's problem in finite
fields with Dickson polynomials. Theory and applications of finite
fields, 185--192, Contemp. Math., 518, Amer. Math. Soc.,
Providence, RI, 2010.



\bibitem{GW} D. G\'{o}mez-P\'{e}rez and A. Winterhof. Multiplicative character sums of
Fermat quotients and pseudorandom sequences, Period. Math. Hungar.
64 (2012) 161--168.

\bibitem{H-B} R. Heath-Brown. An estimate for Heilbronn's exponential sum. Analytic Number Theory:
Proc. Conf.  in honor of Heini Halberstam, Birkh{\"a}user, Boston,
1996, 451--463.

\bibitem{LW} Y. R. Liu and T. D. Wooley. Waring's problem in function fields. J.
Reine Angew. Math. 638 (2010)  1--67.



\bibitem{OS} A. Ostafe  and  I. E. Shparlinski. Pseudorandomness and dynamics of
Fermat quotients. SIAM J. Discr. Math. 25 (2011) 50--71.

\bibitem{OS11}
A. Ostafe  and  I. E. Shparlinski. On the Waring problem with
Dickson polynomials in finite fields. Proc. Amer. Math. Soc. 139
(2011) 3815--3820.


\bibitem{OTW}
A. Ostafe, D. Thomson and  A Winterhof. On the Waring problem with
multivariate Dickson polynomials. Theory and applications of finite
fields, 153--161, Contemp. Math., 579, Amer. Math. Soc., Providence,
RI, 2012.


\bibitem{SS} J. Sauerberg and L. Shu. Fermat quotients over function
fields. Finite Fields Appl. 3 (1997) 275--286.


\bibitem{Sha}M. Sha. The arithmetic of Carmichael quotients. http://arxiv.org/arXiv:1108.2579, 2011.

\bibitem{Shk}I. D. Shkredov. On Heilbronn's exponential sum. Quart. J. Math. 64 (2013) 1221--1230.





\bibitem{S}I. E. Shparlinski. Character sums with Fermat quotients. Quart. J. Math. 62 (2011) 1031--1043.

\bibitem{S2010} I. E. Shparlinski.  Bounds of multiplicative character sums with
Fermat quotients of primes. Bull. Aust. Math. Soc. 83 (2011)
456--462.

\bibitem{S2011} I. E. Shparlinski.  On the value set of Fermat
quotients. Proc. Amer. Math. Soc. 140 (2012) 1199--1206.

\bibitem{S2011b} I. E. Shparlinski.  Fermat quotients: Exponential sums, value set and primitive
roots. Bull. Lond. Math. Soc.  43 (2011) 1228--1238.

\bibitem{SW}
I. E. Shparlinski  and A. Winterhof. Distribution of values of
polynomial Fermat quotients. Finite Fields Appl. 19 (2013)
93--104.



\bibitem{TV} T. Tao and V. Vu. Additive Combinatorics. Cambridge Univ. Press,
Cambridge, 2006.



\bibitem{W98} A. Winterhof. On Waring's problem in finite fields. Acta Arith. 87
(1998)  171--177.


\bibitem{W01}
A. Winterhof. A note on Waring's problem in finite fields. Acta
Arith. 96 (2001) 365--368.


\bibitem{WW}
A. Winterhof and  C.  van de Woestijne. Exact solutions to Waring's
problem for finite fields. Acta Arith. 141 (2010) 171--190.






\end{thebibliography}
\end{document}